\documentclass[twocolumn,twoside]{IEEEtran}

\usepackage{amsmath,amssymb}
\usepackage{graphicx,color,psfrag}
\usepackage{palatino}

\newtheorem{theorem}{Theorem}[section]
\newtheorem{lemma}[theorem]{Lemma}
\newtheorem{remark}[theorem]{Remark}
\newtheorem{proposition}[theorem]{Proposition}

\newcommand\real{{\ensuremath{\mathbb{R}}}}

\newcommand{\sphere}{{\mathbb{S}}}

\newcommand{\scirc}{\circ}    

\newcommand{\mathscr}[1]{\boldsymbol{\EuScript{#1}} }

\newcommand{\spn}{\operatorname{span}}
\newcommand\map[3]{#1\colon#2\rightarrow#3}

\newcommand{\invclos}[1]{\operatorname{\overline{Lie}}(#1)}

\newcommand{\atantwo}[2]{\operatorname{arctan2}\left(#1,#2\right)}

\newcommand\mapto[2]{#1\longmapsto#2}

\renewcommand{\t}{\theta}
\newcommand{\f}{\text{f}}

\newcommand\lie[1]{{\ensuremath{\mathfrak{#1}}}}
\newcommand{\so}[1]{\lie{so}(#1)}
\newcommand{\se}[1]{\lie{se}(#1)}
\newcommand{\Lie}[1]{{\ensuremath{\mathrm{#1}}}}
\newcommand{\SO}[1]{{\ensuremath{\mathrm{SO}(#1)}}}
\newcommand{\SE}[1]{{\ensuremath{\mathrm{SE}(#1)}}}

\newcommand{\Id}{\operatorname{id}}
\newcommand{\Sign}{\operatorname{sign}}
\newcommand{\Ind}{{\operatorname{ind}}}
\newcommand{\forwardkin}{{\mathcal{FK}}}
\newcommand{\inversekin}{{\mathcal{IK}}}

\newcommand\setdef[2]{\{#1 \; |\enspace#2\}}

\renewcommand{\theenumi}{(\roman{enumi})}

\title{A catalog of inverse-kinematics planners\\ for underactuated systems
  on matrix groups}

\author{Sonia Mart{\'\i}nez\thanks{Departamento de Matem\'atica
    Aplicada IV, Universidad Polit\'ecnica de Catalu\~na, Av. V.
    Balaguer, s/n, Vilanova i la Geltr\'u, 08800, Spain,
    email:~{\texttt{soniam@mat.upc.es}}}, Jorge Cort\'es, and
  Francesco Bullo\thanks{Coordinated Science Laboratory, University of
    Illinois at Urbana-Champaign, Urbana, IL 61801, United States,
    email:~{\texttt{\{jcortes,bullo\}@uiuc.edu}}}}

\begin{document}
\maketitle
\begin{abstract}
  This paper presents motion planning algorithms for underactuated systems
  evolving on rigid rotation and displacement groups.  Motion planning is
  transcribed into (low-dimensional) combinatorial selection and
  inverse-kinematics problems.  We present a catalog of solutions for all
  underactuated systems on $\SE{2}$, $\SO{3}$ and $\SE{2}\times\real$
  classified according to their controllability properties.
\end{abstract}

\section{Introduction}
This paper presents motion planning algorithms for kinematic models of
underactuated mechanical systems; we consider kinematic (i.e.,
driftless) models that are switched control systems, that is,
dynamical systems described by a family of admissible vector fields
and a control strategy that governs the switching between them.  
In particular, we focus on families of left-invariant vector fields
defined on rigid displacements subgroups.

This class of systems arises in the context of kinematic modeling and
kinematic reductions for mechanical control systems; see the recent
works~\cite{HA-KT-NS:98,ADL:98c,KML-NS-HA-KT:00,FB-KML:01a,FB-ADL-KML:02a}.
Clearly, the transcription into kinematic models simplifies the motion
planning problem; e.g., \cite{FB-KML:01a} discusses 3R planar manipulators
and~\cite{FB-ADL:01t,SI-KML:03a} discuss the snakeboard system.

\subsection*{Literature review}
Motion planning for kinematic models, sometimes referred to as
driftless or nonholonomic models, is a classic problem in robotics;
see~\cite{DH-LEK-JCL:97a} and also the references therein. In
particular, the algorithms
in~\cite{GL-HJS:93,NEL-PSK:95,RMM-ZXL-SSS:94} focus on dynamical
aspects and exploit controllability properties.

For the class of systems of interest in this paper, the search for a motion
planning algorithm is closely related to the inverse-kinematics problem.
Example inverse-kinematics methods include (i) iterative numerical methods
for nonlinear optimization, see~\cite{AG-BB-RF:85}, (ii) geometric and
decoupling methods for classes of manipulators,
see~\cite{MWS-MV:89,LS-BS:00}, (iii) the Paden-Kahan subproblems approach,
see~\cite{BP:86,RMM-ZXL-SSS:94}, and (iv) the general polynomial
programming approach, see~\cite{DM-JFC:94}.  The latter and more general
method is based on tools from algebraic geometry and relies on
simultaneously solving systems of algebraic equations.  Despite these
efforts, no general methodology is currently available to solve these
problems in closed-form. Accordingly, it is common to provide and catalog
closed-form solutions for classes of relevant example systems;
see~\cite{RMM-ZXL-SSS:94,MWS-MV:89,LS-BS:00}.
  
\subsection*{Problem statement} 
We consider left-invariant control systems evolving on a matrix Lie
subgroup $\Lie{G}\subset\SE{3}$.  Examples include systems on
$\SE{2}$, $\SO{3}$ and $\SE{2} \times \real$.  As usual in Lie group
theory, we identify left-invariant vector fields with their value at
the identity.  Given a family of left-invariant vector fields
$\{V_1,\dots,V_m\}$ on $\Lie{G}$, consider the associated driftless
control system
\begin{equation} \label{eq:kinematicsystem} 
  \dot{g}(t) = \sum_{i=1}^m V_i(g(t)) u_i(t) \, ,
\end{equation}
where $g:\real\rightarrow\Lie{G}$ and where the controls
$(u_1,\dots,u_m)$ take value in $\{(\pm 1,0,\ldots,0),(0,\pm
1,0,\ldots,0),\ldots, (0,\ldots,0,\pm 1)\}$.  For these systems,
controllability can be assessed by algebraic means: it suffices to
check the lack of involutivity of $\spn\{V_1,\dots,V_m\}$.  Recall
that for matrix Lie algebras, Lie brackets are matrix commutators
$[A,B] = AB - BA$.
  
This paper addresses the problem of how to compute feasible motion
plans for the control system~\eqref{eq:kinematicsystem} by
concatenating a finite number of flows along the input vector fields.
We call a flow along any input vector field a \emph{motion primitive} and its
duration a \emph{coasting time}.  Therefore, motion planning is
reduced to the problem of selecting a finite-length combination of $k$
motion primitives $(i_1,\dots,i_k) \in \{1,\dots,m\}^k$\ and computing
appropriate coasting times $(t_1,\dots,t_k)\in\real^k$\ that steer the
system from the identity in the group to any target configuration
$g_{\f}\in\Lie{G}$. In mathematical terms, we need to solve
\begin{equation*}
  g_{\f} = \exp(t_1 V_{i_1}) \cdots \exp(t_k V_{i_k}) .
\end{equation*}
Hence, motion planning is transcribed into low-dimensional combinatorial
selection and inverse-kinematics problems.

\subsection*{Contribution}

The contribution of this paper is a catalog of solutions for
underactuated example systems defined on $\SE{2}$, $\SO{3}$, or
$\SE{2}\times\real$.  Based on a controllability analysis, we classify
families of underactuated systems that pose qualitatively different
planning problems.  For each family, we solve the planning problem by
providing \emph{a} combination of $k$ motion primitives and corresponding
closed-form expressions for the coasting times.  In each case, we
attempt to select $k=\dim(\Lie{G})$: generically, this is the minimum
necessary (but sometimes not sufficient) number of motion primitives needed.
If the motion planning algorithm entails exactly $\dim(\Lie{G})$
motion primitives, i.e., minimizes the number of switches, we will refer to it
as a \emph{switch-optimal} algorithm.
Sections~\ref{sec:SE2},~\ref{sec:SO3}, and~\ref{sec:SE2timesR} present
switch-optimal planners for $\SE{2}$, $\SO{3}$, and
$\SE{2}\times\real$, respectively.

\subsection*{Notation}

Here we briefly collect the notation used throughout the paper. Let
$S$ be a set, $\map{\Id_{S}}{S}{S}$ denote the identity map on $S$ and
let $\map{\Ind_S}{\real}{\real}$ denote the characteristic function of
$S$, i.e., $\Ind_S(x)=1$ if $x \in S$ and $\Ind_S (x)=0$ if $x \not
\in S$.  Let $\atantwo{x}{y}$ denote the arctangent of ${y}/{x}$
taking into account which quadrant the point $(x,y)$ is in. We make
the convention $\atantwo{0}{0}=0$.  Let $\map{\Sign}{\real}{\real}$ be
the sign function, i.e., $\Sign (x)=1$ if $x>0$, $\Sign (x)=-1$ if
$x<0$ and $\Sign (0)=0$.  Let $A_{ij}$ be the $(i,j)$ element of the
matrix $A$. Given $v,w\in\real^n$, let $\arg(v,w)\in[0,\pi]$ denote
the angle between them.  Let $\|\cdot\|$ denote the Euclidean norm.

Given a family of left-invariant vector fields $\{V_1,\dots,V_m\}$ on
$\Lie{G}$, we associate to each multiindex
$(i_1,\dots,i_k)\in\{1,\dots,m\}^k$ the forward-kinematics map
$\forwardkin^{(i_1,\dots,i_k)}: \real^k \rightarrow \Lie{G}$ given by
$(t_1,\dots,t_k) \mapsto \exp(t_1 V_{i_1}) \cdots \exp(t_k V_{i_k})$.

\section{Catalog for $\SE{2}$}\label{sec:SE2}
Let $\left\{ e_{\t}, e_x, e_y \right\}$ be the basis of $\se{2}$:
\begin{equation*} 
  e_{\t} =  \begin{bmatrix} 
      0 & -1 & 0 \\
      1 & 0 & 0 \\
      0 & 0 & 0
    \end{bmatrix}, \,
  e_x = \begin{bmatrix}
    0 & 0 & 1 \\
    0 & 0 & 0 \\
    0 & 0 & 0
  \end{bmatrix},  \,
  e_y = \begin{bmatrix} 
    0 & 0 & 0 \\
    0 & 0 & 1 \\
    0 & 0 & 0
  \end{bmatrix} 
  .
\end{equation*} 
Then, $[e_{\t},e_x]=e_y$, $[e_y,e_{\t}] = e_x$ and $[e_x,e_y]=0$.  For ease
of presentation, we write $V\in\se{2}$ as $ V = a e_\t + b e_x + c e_y
\equiv (a,b,c)$, and $g\in\SE{2}$ as
\begin{equation*}
  g =  
 \begin{bmatrix}
   \cos \t & - \sin \t & x \\
   \sin \t & \cos \t & y \\
   0 & 0 & 1
 \end{bmatrix}
 \equiv (\t,x,y)\,.
\end{equation*} 
With this notation, $\map{\exp}{\se{2}}{\SE{2}}$ is
\begin{multline*}
  \exp(a,b,c) \\
  = \left(a \, , \; \frac{\sin a}{a}b - \frac{1 - \cos a}{a}c\,,\; \frac{1 -
      \cos a}{a}b + \frac{\sin a}{a} c\right)
\end{multline*}
for $a \neq 0$, and $\exp(0,b,c)=(0,b,c)$.

\begin{lemma}(Controllability conditions). \label{lem:ctrl-SE2}
  Consider two left-invariant vector fields $V_1 = (a_1,b_1, c_1)$ and $V_2
  = (a_2,b_2, c_2)$ in $\se{2}$.  Their Lie closure is full rank if and
  only if $a_1 b_2 -b_1 a_2 \neq 0$ or $c_1 a_2 - a_1 c_2 \neq 0$.
\end{lemma}

\begin{proof}
  Given $[V_1,V_2] = (0,c_1 a_2 - a_1 c_2 \,,\, a_1 b_2 - b_1 a_2)$, one can
  see that $\spn{\left\{V_1, V_2, [V_1,V_2]\right\}} = \se{2}$ if and only~if
  \begin{multline*}
    \det
    \begin{bmatrix}
      a_1 & b_1 & c_1  \\
      a_2 & b_2 & c_2 \\
      0 & c_1 a_2 - c_2a_1 & b_2 a_1- b_1 a_2 
    \end{bmatrix} \\ 
    = ( a_1 b_2 -b_1 a_2 )^2 + ( c_1 a_2 - a_1 c_2)^2 \neq 0 \,.
  \end{multline*}
\end{proof}

Let $V_1 = (a_1,b_1, c_1)$ and $V_2 = (a_2,b_2, c_2)$ satisfy the
controllability condition in Lemma~\ref{lem:ctrl-SE2}. Accordingly, either
$a_1$ or $a_2$ is different from zero. Without loss of generality, we will
assume that $a_1 \neq 0$, and take $a_1=1$. As a consequence of
Lemma~\ref{lem:ctrl-SE2}, there are two qualitatively different cases to be
considered:
\begin{align*}
  {\mathcal{S}}_1 & = \{(V_1,V_2) \in \se{2} \times \se{2} \, | \;
  V_1=(1,b_1,c_1),\\
  & \hspace*{3.5cm} V_2 =(0,b_2,c_2) \, \text{and} \, b_2^2 + c_2^2 = 1 \}
  \,,   \\
  {\mathcal{S}}_2 &= \{(V_1,V_2) \in \se{2} \times \se{2} \, | \; V_1=(1,b_1,
  c_1), \\
  & \hspace*{1.6cm} V_2 =(1,b_2, c_2) \, \text{and either} \, b_1 \neq b_2 \,
  \text{or} \, c_1 \neq c_2 \} \, .
\end{align*}
Since $\dim(\se{2})=3$, we need at least three motion primitives along the flows of
$\left\{V_1,V_2\right\}$ to plan any motion between two desired
configurations. Consider the map
$\map{\forwardkin^{(1,2,1)}}{\real^3}{\SE{2}}$.
In the following propositions, we compute solutions for ${\mathcal{S}}_1$ and
${\mathcal{S}}_2$-systems.

\begin{proposition} (Inversion for ${\mathcal{S}}_1$-systems on $\SE{2}$). 
  \label{prop:SE(2)-case-A}
  Let $(V_1,V_2) \in {\mathcal{S}}_1$.  Consider the map
  $\map{\inversekin[{\mathcal{S}}_1]}{\SE{2}}{\real^3}$,
  \begin{equation*}
    \inversekin[{\mathcal{S}}_1] (\theta,x,y) 
    = (\atantwo{\alpha}{\beta}, \rho, \theta -\atantwo{\alpha}{\beta}), 
  \end{equation*}
  where $\rho=\sqrt{\alpha^2 + \beta^2}$ and 
  \begin{align*}
    \begin{bmatrix}
      \alpha\\
      \beta
    \end{bmatrix}
    &=
    \begin{bmatrix}
      b_2 & c_2\\
      -c_2 & b_2
    \end{bmatrix}
    \left(
    \begin{bmatrix}
      x \\ y
    \end{bmatrix}
    -
    \begin{bmatrix}
      -c_1 & b_1 \\ 
      b_1 & c_1     
    \end{bmatrix}
    \begin{bmatrix}
      1 - \cos \t \\
      \sin \t
    \end{bmatrix}    
  \right).
  \end{align*}
  Then, $\inversekin[{\mathcal{S}}_1]$ is a global right inverse of
  $\forwardkin^{(1,2,1)}$, that is, it satisfies $\map{\forwardkin^{(1,2,1)}
    \scirc \inversekin[{\mathcal{S}}_1] = \Id_{\SE{2}}}{\SE{2}}{\SE{2}}$.
\end{proposition}
Note that the algorithm provided in the proposition is not only
switch-optimal, but also works globally.

\begin{proof}
  The proof follows from the expression of the map
  $\forwardkin^{(1,2,1)}$.  Let $\forwardkin^{(1,2,1)} (t_1,t_2,t_3) =
  (\theta,x,y)$,
  \begin{align*} 
    \t &= t_1 + t_3 \,, \\
    \begin{bmatrix}
      x \\ y 
    \end{bmatrix}
    & = 
    \begin{bmatrix}
      -c_1 & b_1 \\ 
      b_1 & c_1     
    \end{bmatrix}
    \begin{bmatrix}
      1 - \cos \t \\
      \sin \t
    \end{bmatrix}    
    +  
    \begin{bmatrix}
      b_2 & -c_2 \\
      c_2 & b_2
    \end{bmatrix}
    \begin{bmatrix}
      \cos t_1 \\ 
      \sin t_1     
    \end{bmatrix} t_2\,.  
  \end{align*} 
  The equation in $[x,y]^T$ can be rewritten as $[\alpha, \beta]^T = [ \cos
  t_1, \sin t_1]^T t_2$. The selection $t_1 =\atantwo{\alpha}{\beta}$, $t_2 =
  \rho$ solves this equation.
\end{proof}

\begin{proposition} (Inversion for ${\mathcal{S}}_2$-systems on $\SE{2}$). 
  \label{prop:SE(2)-case-B}
  Let $(V_1,V_2) \in {\mathcal{S}}_2$. Define the neighborhood of the
  identity in~$\SE{2}$
  \begin{multline*}    
  U = \{(\theta,x,y)\in \SE{2}\,|\enspace \| ( c_1 - c_2, b_1 - b_2 )\|^2\geq 
  \\ 
    \max\{ \| (x ,y ) \|^2 \,, 2(1 - \cos \theta) \|  (b_1,c_1)\|^2 \}. 
  \end{multline*}
  Consider the map
  $\map{\inversekin[{\mathcal{S}}_2]}{U\subset\SE{2}}{\real^3}$ whose
  components are
  \begin{align*}
    \inversekin[{\mathcal{S}}_2]_1 (\theta,x,y) &= \atantwo{\rho}{\sqrt{4
        -\rho^2}} +
    \atantwo{\alpha}{\beta} \,, \\
    \inversekin[{\mathcal{S}}_2]_2 (\theta,x,y) &= \atantwo{2-\rho^2}{\rho
      \sqrt{4-\rho^2}}
    \!, \\
    \inversekin[{\mathcal{S}}_2]_3 (\theta,x,y) &= \theta
    -\inversekin[{\mathcal{S}}_2]_1
    (\theta,x,y)-\inversekin[{\mathcal{S}}_2]_2 (\theta,x,y) \, ,
  \end{align*}
  and $\rho=\sqrt{\alpha^2 + \beta^2}$ and 
  \begin{align*}
    \begin{bmatrix}
      \alpha\\
      \beta
    \end{bmatrix}
    &= \frac{1}{ \| ( c_1 - c_2, b_1 - b_2 )\|^2}
    \begin{bmatrix}
      c_1 - c_2 & b_2 - b_1\\
      b_1 - b_2 & c_1 - c_2
    \end{bmatrix}
    \\
    & \hspace*{2cm} \cdot \left(
    \begin{bmatrix}
      x \\ y
    \end{bmatrix}
    -
    \begin{bmatrix}
      -c_1 & b_1 \\
      b_1 & c_1
    \end{bmatrix}
    \begin{bmatrix}
      1 - \cos \t \\
      \sin \t
    \end{bmatrix}    
  \right) \, .
  \end{align*}
  Then, $\inversekin[{\mathcal{S}}_2]$ is a local right inverse of
  $\forwardkin^{(1,2,1)}$, that is, it satisfies $\map{\forwardkin^{(1,2,1)}
    \scirc \inversekin[{\mathcal{S}}_2] = \Id_{U}}{U}{U}$.
\end{proposition}

\begin{proof}
  If $(\theta,x,y) \in U$, then
  \begin{multline*}
    \rho = \| (\alpha,\beta)\| \le \frac{1}{\| ( c_1 - c_2, b_1 - b_2 )\|
      } \\
    \cdot \left( \|(x, y)\| + \Big \|
      \begin{bmatrix}
        -c_1 & b_1 \\
        b_1 & c_1
      \end{bmatrix}
      \begin{bmatrix}
        1 - \cos \t  \\
        \sin \t 
      \end{bmatrix}
      \Big \| \right) \leq 2 \, ,
  \end{multline*}
  and hence $\inversekin[{\mathcal{S}}_2]$ is well-defined on $U$.  Let
  $\inversekin[{\mathcal{S}}_2] (\theta,x,y)=(t_1,t_2,t_3)$. The components
  of $\forwardkin^{(1,2,1)} (t_1,t_2,t_3)$ are
  \begin{align*}
    & \forwardkin^{(1,2,1)}_1 (t_1,t_2,t_3) = t_1 + t_2 + t_3 \,, \\
    & \begin{bmatrix}
      \forwardkin^{(1,2,1)}_2 (t_1,t_2,t_3) \\
      \forwardkin^{(1,2,1)}_3 (t_1,t_2,t_3)
    \end{bmatrix}
    =
    \begin{bmatrix}
      -c_1 & b_1 \\
      b_1 & c_1
    \end{bmatrix}
    \begin{bmatrix}
      1 - \cos \t  \\
      \sin \t
    \end{bmatrix}
    \\
    & \hspace*{1.5cm} +
    \begin{bmatrix}
      c_1-c_2 &  b_1-b_2  \\
      b_2- b_1 & c_1 - c_2  
    \end{bmatrix}
    \begin{bmatrix}
      \cos t_1 - \cos (t_1 + t_2)\\
      \sin t_1 - \sin (t_1 + t_2)
    \end{bmatrix} \,.
  \end{align*}
  In an analogous way to the previous proof, one verifies
  $\forwardkin^{(1,2,1)} (t_1,t_2,t_3) = (\theta,x,y)$.
\end{proof}

\begin{remark}{\rm
    The map $\inversekin[{\mathcal{S}}_2]$ in
    Proposition~\ref{prop:SE(2)-case-B} is a local right inverse to
    $\forwardkin^{(1,2,1)}$ on a domain that strictly contains $U$. In other
    words, our estimate of the domain of $\inversekin[{\mathcal{S}}_2]$ is
    conservative.  For instance, for points of the form $(0,x,y) \in \SE{2}$,
    it suffices to ask for
    \begin{equation*}
      \| (x,y) \| \le 2 \| ( c_1 - c_2, b_1 - b_2 )\| \,.
    \end{equation*}
    For a point $(\theta,0,0) \in \SE{2}$, it suffices to ask for
    \begin{equation*}
      (1 - \cos\theta) \|(b_1,c_1)\|^2 \le 2 \| ( c_1 - c_2, b_1 - b_2 )\|^2
      \,. 
    \end{equation*}    
    Additionally, without loss of generality, it is convenient to assume
    that the vector fields $V_1$, $V_2$ satisfy $b_1^2 + c_1^2 \le b_2^2 +
    c_2^2$, so as to maximize the domain $U$. }
\end{remark}

\begin{figure}[htbp]
  \centering \includegraphics[width=.4\linewidth]{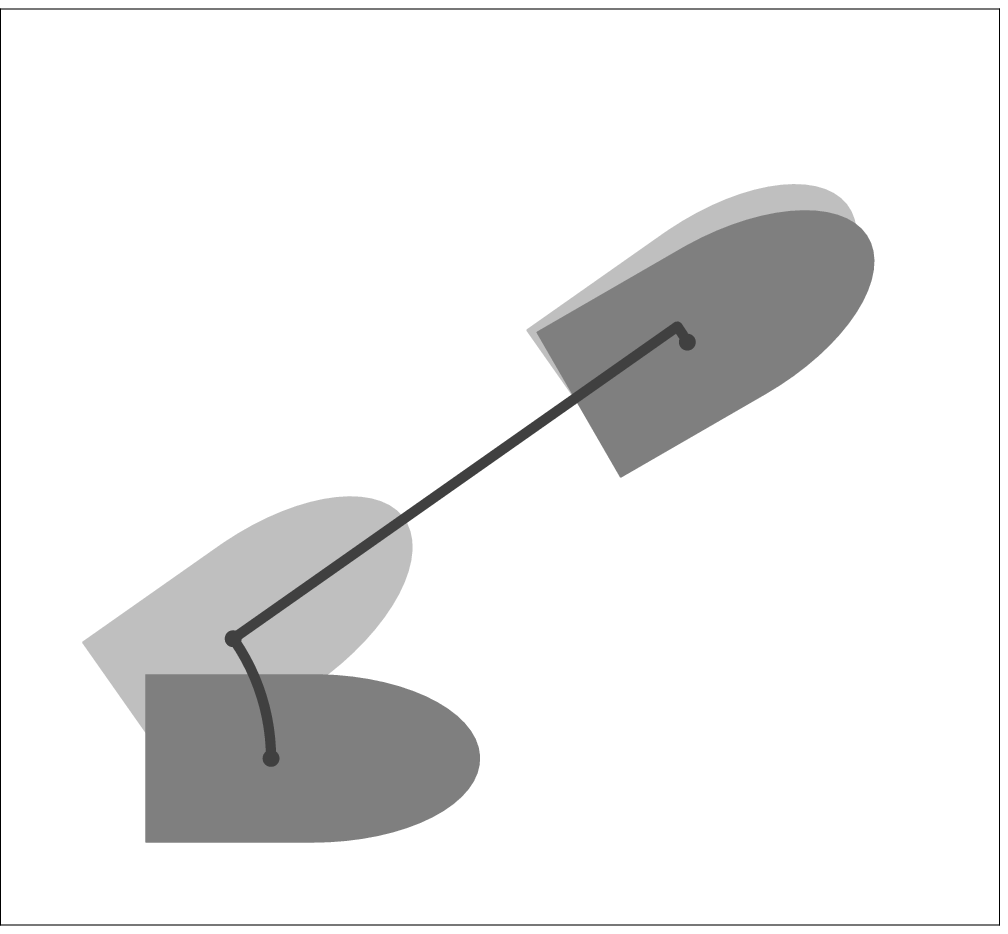} \quad
  \includegraphics[width=.4\linewidth]{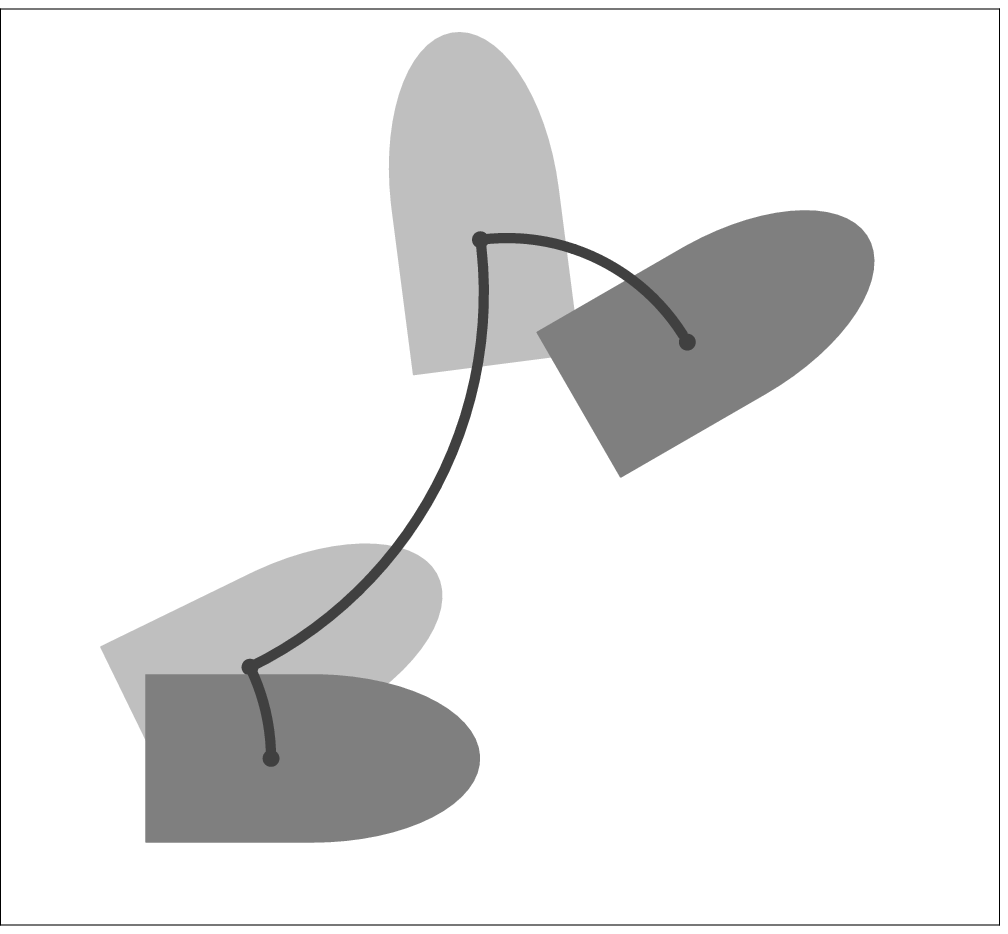}
  \caption{We illustrate the inverse-kinematics planners for
    ${\mathcal{S}}_1$ and ${\mathcal{S}}_2$-systems.  The parameters of
    both systems are $(b_1,c_1)=(0,.5)$, $(b_2,c_2)=(1,0)$. The target
    location is $(\pi/6,1,1)$. Initial and target location are depicted in
    dark gray.}
  \label{fig:se2}
\end{figure}

We illustrate the performance of the algorithms in Fig.~\ref{fig:se2}.

\section{Catalog for SO(3)}\label{sec:SO3}
Let $\left\{ \widehat{e}_x, \widehat{e}_y, \widehat{e}_z \right\}$ be the
basis of $\so{3}$:
\begin{gather*}
  \widehat{e}_x =\!
  \begin{bmatrix} 
    0 & 0 & 0 \\
    0 & 0 & -1 \\
    0 & 1 & 0
  \end{bmatrix} 
  \hspace*{-3pt} , \, \widehat{e}_y =\!
  \begin{bmatrix} 
    0 & 0 & 1 \\
    0 & 0 & 0 \\
    -1 & 0 & 0
  \end{bmatrix} 
  \hspace*{-3pt} , \, \widehat{e}_{z} =\!
  \begin{bmatrix} 
    0 & -1 & 0 \\
    1 & 0 & 0 \\
    0 & 0 & 0
  \end{bmatrix} 
  \hspace*{-3pt} .
\end{gather*} 
Here we make use of the notation $\widehat{V}=a \widehat{e}_x +
b\widehat{e}_y + c \widehat{e}_z \equiv \widehat{(a, b, c)}$ based on the
Lie algebra isomorphism $\widehat{\cdot}: (\real^3, \times ) \rightarrow
(\so{3}, [\cdot,\cdot])$.  Rodrigues formula~\cite{RMM-ZXL-SSS:94} for the
exponential $\exp : \so{3} \rightarrow \SO{3}$ is
\begin{equation*}
  \exp  (\widehat{\eta}) =  I_3  +  \frac{\sin \|  \eta  \|}{\| \eta  \|}
  \widehat{\eta}    +   \frac{1    -   \cos    \|   \eta    \|}{\|   \eta
    \|^2}\widehat{\eta}^2\,.
\end{equation*} 
The commutator relations are $\left[ \widehat{e}_x,\widehat{e}_{z} \right] =
-\widehat{e}_y$, $\left[ \widehat{e}_y,\widehat{e}_{z} \right] =
\widehat{e}_x$ and $\left[ \widehat{e}_x,\widehat{e}_{y} \right] =
\widehat{e}_z$.

\begin{lemma} (Controllability conditions). \label{lem:ctrl-SO3}
  Consider two left-invariant vector fields $V_1 = (a_1, b_1, c_1)$ and $V_2
  = (a_2, b_2, c_2)$ in $\so{3}$.  Their Lie closure is full rank if and only
  if $c_1 a_2 -a_1 c_2 \neq 0$ or $b_1 c_2 - c_1 b_2 \neq 0$ or $b_1 a_2 -
  a_1 b_2 \neq 0$.
\end{lemma}
\begin{proof}
  Given the equality $[\widehat{V}_1,\widehat{V}_2] = \widehat{V_1 \times
    V_2}$, with $V_1 \times V_2 = ( b_1 c_2-b_2 c_1, c_1 a_2 -c_2 a_1, a_1
  b_2 - a_2 b_1)$, one can see that $\spn{\left\{V_1, V_2, [V_1,V_2]\right\}}
  = \so{3}$ if and only if
  \begin{multline*}
    \det
    \begin{bmatrix}
      a_1 & b_1 & c_1 \\
      a_2 & b_2 & c_2 \\
      b_1 c_2 - b_2 c_1 & c_1 a_2- c_2 a_1 & a_1 b_2- a_2 b_1
    \end{bmatrix}
    = \\
    ( b_1 c_2 - b_2 c_1)^2 + ( c_1 a_2 - c_2 a_1)^2 + ( a_1 b_2 - a_2 b_1)^2
    \neq 0 \,.
  \end{multline*}
\end{proof}

Let $V_1$, $V_2$ satisfy the controllability condition in
Lemma~\ref{lem:ctrl-SO3}. Without loss of generality, we can assume
$V_1={e}_z$ (otherwise we perform a suitable change of coordinates), and
$\|V_2 \|=1$. In what follows, we let $V_2=(a,b,c)$. Since ${e}_z$ and
$V_2$ are linearly independent, necessarily $a^2 + b^2 \neq 0 $ and $c\neq\pm
1$. Since $\dim(\so{3})=3$, we need at least three motion primitives to plan any
motion between two desired configurations.  Consider the map
$\map{\forwardkin^{(1,2,1)}}{\real^3}{\SO{3}}$, that is
\begin{equation}\label{eq:exp-so3}
  \forwardkin^{(1,2,1)}(t_1,t_2,t_3) = \exp (t_1 \widehat{e}_z) \exp(t_2
  \widehat{V}_2)  \exp  (t_3 \widehat{e}_z) \, .
\end{equation}
Observe that equation~\eqref{eq:exp-so3} is similar to the formula for
certain sets of Euler angles; see~\cite{RMM-ZXL-SSS:94}.

\begin{proposition} (Inversion for systems on $\SO{3}$).
  Let $V_1=(0, 0, 1)$ and $V_2 =(a,b,c)$, with $a^2+b^2\neq 0 $ and $c \neq
  \pm 1$. Define the neighborhood of the identity in $\SO{3}$
  \begin{align*}    
    U = \setdef{R \in\SO{3}}{R_{33} \in [2c^2 -1,1]}.
  \end{align*}
  Consider the map $\map{\inversekin}{U \subset \SO{3}}{\real^3}$ whose
  components are
  \begin{align*}
    \inversekin_1 (R) &= \atantwo{w_1 R_{13} + w_2 R_{23}}{ -w_2 R_{13} + w_1
      R_{23}} \,,\\
    \inversekin_2 (R) &= \arccos\left(\frac{R_{33} - c^2}{1 - c^2}\right) \,,
    \\
    \inversekin_3 (R) &= \atantwo{v_1 R_{31} + v_2 R_{32}}{v_2 R_{31} - v_1
      R_{32}} \, ,
  \end{align*}
  where, for $z=(1-\cos(\inversekin_2 (R)), \sin(\inversekin_2 (R)))^T$,
  \begin{align*}
    \begin{bmatrix}
      w_1\\
      w_2
    \end{bmatrix}
    &=
    \begin{bmatrix}
      a c & b\\
      c b & - a
    \end{bmatrix}  z 
    \,, \quad
    \begin{bmatrix}
      v_1\\
      v_2
    \end{bmatrix}    =
    \begin{bmatrix}
      a c & - b\\
      c b & a
    \end{bmatrix} z.
  \end{align*}
  Then, $\inversekin$ is a local right inverse of $\forwardkin^{(1,2,1)}$,
  that is, it satisfies $\forwardkin^{(1,2,1)} \scirc \inversekin = \Id_{U}:U
  \to U$.
\end{proposition}

\begin{proof}
  Let $R \in U$. Then, $|\frac{R_{33} - c^2}{1 - c^2} | \le 1$, and
  hence $\inversekin(R)$ is well-defined. Denote
  $t_i=\inversekin_i(R)$ and let us show
  $R=\forwardkin^{(1,2,1)}(t_1,t_2,t_3)$.  Recall that the rows
  (resp.~the columns) of a rotation matrix consist of orthonormal
  vectors in $\real^3$. Therefore, the matrix
  $\forwardkin^{(1,2,1)}(t_1,t_2,t_3) \in \SO{3}$ is determined by its
  third column $\forwardkin^{(1,2,1)}(t_1,t_2,t_3) e_{z}$ and its
  third row $e_z ^T \forwardkin^{(1,2,1)}(t_1,t_2,t_3)$.  The factors
  in~\eqref{eq:exp-so3} admit the following closed-form expressions.
  For $c_t=\cos t$ and $s_t=\sin t$,
  \begin{equation*}
    \exp(t  \widehat{e}_z) =
    \begin{bmatrix}
      c_t & -s_t & 0 \\
      s_t & c_t & 0 \\
      0 & 0 & 1
    \end{bmatrix} \, ,
  \end{equation*}
  and $\exp(t \widehat{V}_2)$ equals
  \begin{equation*}
    {\small
    \begin{bmatrix} 
      a^2 + (1-a^ 2)c_t& ba(1 - c_t)- c s_t & ca(1 - c_t) +b s_t \\
      ab(1 - c_t) +c s_t & b^2 + (1-b^ 2)c_t & cb(1 - c_t) -a  s_t\\
      ac(1 - c_t) -b s_t & bc(1 - c_t) +a s_t & c^2 + (1-c^2)c_t\\
    \end{bmatrix} \,.}
  \end{equation*}
  Now, using the fact that $\exp (t \widehat{e}_z)e_z = e_z$, we get
  \begin{multline*}
    \forwardkin^{(1,2,1)} (t_1,t_2,t_3) e_z = \exp (t_1 \widehat{e}_z) \exp
    (t_2
    \widehat{V}_2) \exp (t_3 \widehat{e}_z) e_z \\
    = \exp (t_1 \widehat{e}_z) \exp (t_2 \widehat{V}_2) e_z = \exp (t_1
    \widehat{e}_z)
    \begin{bmatrix}
      w_1\\
      w_2\\
      R_{33}
    \end{bmatrix}
    = R e_z \, .
  \end{multline*}
  A similar computation shows that $e_z^T
  \forwardkin^{(1,2,1)}(t_1,t_2,t_3)= e_z^T R$, which concludes the
  proof.
\end{proof}
 
\begin{remark}{\rm
    If $\widehat{e}_z$ and $V_2$ are perpendicular, then $U=SO(3)$ and the
    map $\inversekin$ is a global right inverse of $\forwardkin^{(1,2,1)}$.
    Otherwise, let us provide an equivalent formulation of the constraint
    $R_{33}\in [2 c^2 -1,1]$ in terms of the axis/angle representation of
    the rotation matrix $R$.  Recall that there always exist a, possibly
    non-unique, rotation angle $\t\in[0,\pi]$ and an unit-length axis of
    rotation $\omega\in\sphere^2$ such that $R=\exp(\widehat{\omega}{\t})$.
    Because $\widehat{\omega}^2 = \omega^T \omega - I_3$, an equivalent
    statement of Rodrigues formula is
    \begin{align*}
      R &= I_3 + \widehat{\omega} \sin \t + (1 - \cos \t)(\omega^T \omega -
      I_3).
    \end{align*}
    From $e_z^T\omega = \cos (\arg (e_z, \omega))$, we compute
    \begin{align}
      e_z^T R e_z &= e_z^T e_z + (1 - \cos \t) ((e_z^T \omega)^2 - e^T_z e_z)
      \nonumber
      \\
      &= 1 + (1 - \cos \t) ((e_z^T \omega)^2 - 1)
      \nonumber \\
      &= 1 - \sin^2 (\arg(e_z, \omega))(1 -\cos \t) \, .
      \label{eq:equivalent-U}
    \end{align}
    Therefore, $R_{33}\in [2 c^2 -1,1]$ if and only if
    \begin{multline*}
      1 - \sin^2(\arg(e_z,\omega)) (1 - \cos \t)  \ge 2 c^2 -1 
      \\
      \iff \,
      \sin^2(\arg(e_z,\omega))(1 - \cos \t) \le 2(1 -  c^2)\,.
    \end{multline*}
    Two sufficient conditions are also meaningful. In terms of the rotation
    angle, if $|\t| \le \arccos (2 c^2 -1)$ then $1 - \cos \t \le 2(1 -
    c^2)$, and in turn equation~\eqref{eq:equivalent-U} is satisfied.  In
    terms of the axis of rotation, a sufficient condition for
    equation~\eqref{eq:equivalent-U} is $\sin^2 (\arg(e_z,\omega)) \le
    \sin^2 (\arg(e_z, V_2)) = 1-c^2$. }
\end{remark}

\begin{figure}[htbp]
  \centering
  \includegraphics[width=.9\linewidth]{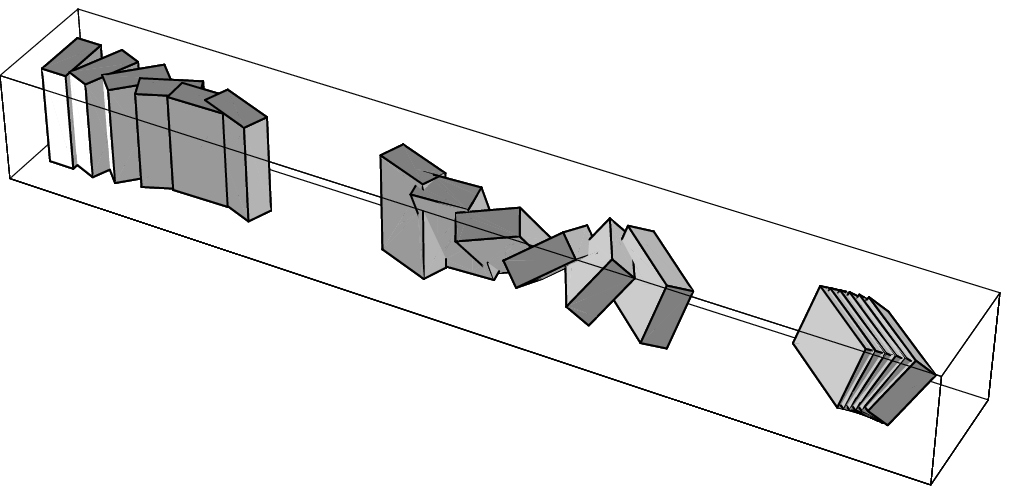}
  \caption{We illustrate the inverse-kinematics planner on $\SO{3}$. 
    The system parameters are $(a,b,c)=(0,1/\sqrt{2},1/\sqrt{2})$. The
    target final rotation is $\exp(\pi/3,\pi/3,0)$. To render the
    sequence of three rotations visible, the body is translated along
    the inertial $x$-axis.}
  \label{fig:so3}
\end{figure}

We illustrate the performance of the algorithm in Fig.~\ref{fig:so3}.

\section{Catalog for $\SE{2} \times \real$}\label{sec:SE2timesR}

Let $\left\{ (e_\theta,0),(e_x,0),(e_y,0),(0,0,0,1) \right\}$ be a
basis of $\se{2} \times \real$, where $\left\{ e_{\theta},e_x,e_y,
\right\}$ stands for the basis of $\se{2}$ introduced in
Section~\ref{sec:SE2}. With a slight abuse of notation, we let
$e_\theta$ denote $(e_\theta,0)$, and we similarly redefine $e_x$ and
$e_y$.  We also let $e_z = (0,0,0,1)$. The only non-vanishing Lie
algebra commutators are $[e_{\t},e_x] = e_y$ and $ [e_{\t},e_y] = -
e_x$.

A left-invariant vector field $V$ in $\se{2} \times \real$ is written
as $V = a e_{\theta} + b e_x + c e_y + d e_z \equiv (a,b,c,d)$, and
$g\in\SE{2} \times \real$ as $g=(\theta,x,y,z)$. The exponential map,
$\exp: \se{2} \times \real \longrightarrow \SE{2} \times \real$, is
given component-wise by the exponential on $\se{2}$ and $\real$,
respectively. That is, $\exp (V)$ is equal to
\begin{equation*}
 \left(a \,,\; \frac{\sin a}{a}b  - \frac{1 - \cos  a}{a}c\,,\;
  \frac{1 - \cos a}{a}b + \frac{\sin a}{a} c\, ,\; d \right)
\end{equation*}
if $a \neq 0$, and $\exp ( V) = (0,b,c,d)$ if $a = 0$.

\begin{lemma}  \label{lem:ctrl-SE2timesR-two}
  (Controllability conditions for systems in $\SE{2}\times\real$ with 2
  inputs). Consider two left-invariant vector fields $V_1 = (a_1,b_1, c_1,
  d_1)$ and $V_2 = (a_2,b_2, c_2, d_2)$ in $\se{2} \times \real$.  Their Lie
  closure is full rank if and only if~$a_2 d_1 - d_2 a_1 \neq 0$, and either
  $c_1 a_2 - a_1 c_2 \neq 0$ or~$a_1 b_2 - b_1 a_2 \neq 0$.
\end{lemma}

\begin{proof}
  Since $[V_1,V_2] = (0\,,c_1 a_2 - a_1 c_2 \,,\, a_1 b_2 - b_1 a_2 \,,
  0)\neq 0$, we deduce that either $c_1 a_2 - a_1 c_2 \neq 0$ or $a_1 b_2 -
  b_1 a_2 \neq 0$. In particular, this implies that necessarily $a_1 \neq 0$
  or $a_2 \neq 0$. Assume $a_1 \neq 0$. Now,
  \begin{align*}
    [V_1,[V_1,V_2]] = (0, a_1(-b_2 a_1 + b_1 a_2) , a_1(c_1 a_2 - c_2 a_1),0)
    \,,
  \end{align*}
  and note that $[V_2,[V_1,V_2]] = (a_2/a_1) [V_1,[V_1,V_2]]$. Finally,
  $\invclos{\{V_1,V_2\}} = \se{2} \times \real$ if and only if
  \begin{multline*}
    \det
    \begin{bmatrix}
      b_1 & c_1 & d_1 & a_1 \\
      b_2 & c_2 & d_2 & a_2 \\
      c_1 a_2 - c_2 a_1 & b_2 a_1 - b_1 a_2 & 0 & 0 \\
      a_1(-b_2 a_1 + b_1 a_2) & a_1(c_1 a_2 - c_2 a_1) & 0 & 0
    \end{bmatrix} =
    \\
    a_1 (a_2 d_1 - d_2 a_1)\left[ (c_1 a_2 - c_2 a_1)^2 + (-b_2 a_1 + b_1 a_2
      )^2 \right] \neq 0.
  \end{multline*}
  Since $[V_1,V_2] \neq 0$, this reduces to $a_2 d_1 - d_2 a_1 \neq 0$.
\end{proof}

Let $V_1$, $V_2$ satisfy the controllability condition in
Lemma~\ref{lem:ctrl-SE2timesR-two}. Without loss of generality, we can assume
$a_1 = 1$. As in the case of $\SE{2}$, there are two qualitatively different
situations to be considered:
\begin{align*}
  {\mathcal{T}}_1 & = \{ (V_1,V_2) \in (\se{2} \times \real)^2 \, | \;
  V_1=(1,b_1, c_1, d_1),\\
  & \hspace*{3.35cm} V_2 =(0,b_2, c_2,1) \, \text{and} \, b_2^2 + c_2^2
  \neq 0 \} \, , \\
  {\mathcal{T}}_2 & = \{ (V_1,V_2) \in (\se{2} \times \real)^2 \, | \;
  V_1=(1,b_1, c_1, d_1), \\
  & V_2 =(1,b_2, c_2, d_2), d_1 \neq d_2 \, \text{and either} \, b_1 \neq b_2
  \, \text{or} \, c_1 \neq c_2\} \, .
\end{align*}

\begin{lemma}\label{lem:ctrl-SE2timesR-three}
  (Controllability conditions for $\SE{2}\times\real$ systems with 3
  inputs).  Consider three left-invariant vector fields
  $V_i=(a_i,b_i,c_i,d_i)$, $i = 1,2,3$ in $\se{2} \times \real$.
  Assume $\invclos{ \left\{ V_{i_1}, V_{i_2}\right\}} \subsetneq
  \se{2} \times \real$, for $i_j \in \left\{1,2,3 \right\}$ and
  $\invclos{ \left\{ V_{1}, V_{2}, V_3 \right\}} = \se{2} \times
  \real$.  Then, possibly after a reordering of the vector fields,
  they must fall in one of the following cases:
  \begin{align*}
    {\mathcal{T}}_3 &= \{ (V_1,V_2,V_3) \in (\se{2} \times \real)^3 \, | \;
    V_1 = (1,b_1,c_1,d_1), V_2 = \\
    & (0,b_2,c_2,0), V_3 = (1,b_1,c_1,d_3), d_1 \neq d_3 \, \text{and} \,
    b_2^2 + c_2^2 \neq 0\} , \\
    {\mathcal{T}}_4 &= \{ (V_1,V_2,V_3) \in (\se{2} \times \real)^3 \, | \;
    V_1 = (1,b_1,c_1,d_1), V_2 = \\
    & (0,b_2, c_2, 0), V_3 \! = \! (0,0,0,d_3), 0 \! \neq \! d_3 \!
    \neq  d_1 \, \text{and} \, b_2^2 + c_2 ^2 \! \neq \! 0\} , \\
    {\mathcal{T}}_5 &= \{ (V_1,V_2,V_3) \in (\se{2} \times \real)^3 \, | \;
    V_1 = (1,b_1,c_1,d_1), \\
    & \hspace*{1.5cm} V_2 = (1,b_2,c_2,d_1), V_3 = (0,0,0,d_3), d_3 \neq 0 \\
    & \hspace*{4cm} \text{and either} \, b_2 \neq b_1 \, \text{or} \, c_1
    \neq c_2\} .
  \end{align*}  
\end{lemma}

\begin{proof}
  Without loss of generality, we can assume that $[V_1,V_2] \neq 0$ and
  $a_1=1$. Since $\invclos{\left\{ V_1,V_2\right\}} \neq \se{2} \times
  \real$, then $a_2 d_1 = d_2$. Given that the Lie closure of $\{
  V_1,V_2,V_3\}$ is full-rank, and $\hbox{dim} (\spn {\left\{
      V_1,V_2,[V_1,V_2]\right\}})= 3$, we have that $d_3 \neq a_3 d_1$. This
  latter fact, together with $\invclos{\left\{V_1,V_3 \right\}} \subsetneq
  \se{2} \times \real$, implies that $[V_1,V_3] = 0$, and therefore $b_3 =
  a_3 b_1$, $c_1 a_3 = c_3$.  We now distinguish two situations depending on
  $[V_2,V_3]$ being zero or not.

  \noindent  \emph{(a) $[V_2,V_3] \neq 0$}.
  Necessarily, $a_3 \neq 0$.  Therefore, we can assume $a_3 = 1$. Since
  $\invclos{\left\{ V_2,V_3 \right\}}$ is not full-rank, then $a_2 =0$.  We
  then have a ${\mathcal{T}}_3 $-system.

  \noindent  \emph{(b) $[V_2,V_3] =0$}.
  Necessarily, $b_3 a_2 = b_2 a_3$ and $c_2 a_3 = c_3 a_2$. Depending on the
  values of $a_2$ and $a_3$, we consider:
    \begin{enumerate}
    \item If $a_2 = a_3 = 0$, then $d_2 = 0$, $d_3 \neq 0$, $b_3 = c_3 = 0$.
      Then, this is a ${\mathcal{T}}_4$-system.
    \item If $a_2 =0$, and $a_3 = 1$, then $b_2 = b_3 a_2 = 0$, $c_2 = c_3a_2
      = 0$ and also $d_2 = d_1 a_2 =0$. This is not possible as it would make
      $V_2 =0$.
    \item If $a_2 = 1$ and $a_3 = 0$, then $b_3=c_3 =0$, and $d_2 = d_1$.
      Therefore, this is a ${\mathcal{T}}_5 $-system.
    \item Finally, if $a_2 =1$ and $a_3 = 1$, then $b_1 = b_2$, $c_1 =c_2$,
      and $d_1 = d_2$, which makes $V_1$ and $V_2$ linearly dependent.
    \end{enumerate}
\end{proof}


\subsection{Two-dimensional input distribution}

Let $V_1$, $V_2$ satisfy the controllability condition in
Lemma~\ref{lem:ctrl-SE2timesR-two}. Since $\dim(\se{2} \times \real)=4$, we
need at least four motion primitives to plan any motion between two desired
configurations.  Consider the map
$\map{\forwardkin^{(2,1,2,1)}}{\real^4}{\SE{2} \times \real}$.

\begin{proposition} \label{prop:SE(2)timesR-case-A-impossibility}
  (Lack of switch-optimal inversion for ${\mathcal{T}}_1$-systems on $\SE{2}
  \times \real$). Let $(V_1,V_2) \in {\mathcal{T}}_1$. Then, the map
  $\forwardkin^{(2,1,2,1)}$ in not invertible at any neighborhood of the
  origin.
\end{proposition}

\begin{proof}
  Let $\forwardkin^{(2,1,2,1)} (t_1,t_2,t_3,t_4) = (\theta,x,y,z)$. Then,
  \begin{align*}
    \theta &= t_2 + t_4\,, \\
    z &= t_1 + t_3 + d_1 (t_2 + t_4) = t_1 + t_3 + d_1 \theta \,,\\
    \begin{bmatrix}
      x \\
      y
    \end{bmatrix}
    &=
    \begin{bmatrix}
      -c_1 \\ b_1
    \end{bmatrix}
    +
    \begin{bmatrix}
      c_1 & b_1 \\
      -b_1 & c_1
    \end{bmatrix}
    \begin{bmatrix}
      \cos \theta \\
      \sin \theta
    \end{bmatrix}
    \\
    & \hspace*{2cm} +
    \begin{bmatrix}
      b_2 \\
      c_2
    \end{bmatrix}
    t_1 +
    \begin{bmatrix}
      b_2 &- c_2 \\
      c_2 & b_2
    \end{bmatrix}
    \begin{bmatrix}
      \cos t_2 \\
      \sin t_2
    \end{bmatrix}
    t_3\,.
  \end{align*}
  Consider a configuration with $\theta =z=0$. Then, the equation in $(x,y)$
  is invertible if and only if the map $f :\real^2 \to  \real^2$ defined by
  \begin{equation*}
    \mapto{
      \begin{bmatrix}
        t_2 \\
        t_3
      \end{bmatrix}
      }{
      \begin{bmatrix}
        \cos t_2 - 1 \\
        \sin t_2
      \end{bmatrix}
      t_3
      }
  \end{equation*}
  is invertible.
  But $f$ can not be inverted in $(0,\beta)$, $\beta \neq 0$.
\end{proof}

\begin{remark}
  {\rm An identical negative result holds if we start taking motion primitives
    along the flow of $V_1$ instead of $V_2$, i.e., if we consider the map
    $\map{\forwardkin^{(1,2,1,2)}}{\real^4}{\SE{2} \times \real}$.}
\end{remark}


Consider the map $\map{\forwardkin^{(1,2,1,2,1)}}{\real^5}{\SE{2} \times
  \real}$.

\begin{proposition} (Inversion for ${\mathcal{T}}_1$-systems on $\SE{2}
  \times \real$).
  \label{prop:SE(2)timesR-case-A}
  Let $(V_1,V_2) \in {\mathcal{T}}_1$. Consider the map
  $\map{\inversekin[{\mathcal{T}}_1]}{\SE{2} \times \real}{\real^5}$ whose
  components are
  \begin{align*}
    \inversekin[{\mathcal{T}}_1]_1 (\theta,x,y,z) &= \pi \Ind_{]-\infty,0[}
    (\gamma-\rho) + \atantwo{\alpha}{\beta}  \\
    & \qquad +  \atantwo{(\rho+\gamma)/2}{0} ,\\
    \inversekin[{\mathcal{T}}_1]_2 (\theta,x,y,z) &= (\gamma-\rho)/2 ,
    \\
    \inversekin[{\mathcal{T}}_1]_3 (\theta,x,y,z) &=
    \atantwo{(\rho^2-\gamma^2)/4}{0} \\
    & + \pi \left( \Ind_{]-\infty,0[} (\gamma+\rho) -\Ind_{]-\infty,0[}
      (\gamma-\rho) \right) , \\
    \inversekin[{\mathcal{T}}_1]_4 (\theta,x,y,z) &= (\gamma+\rho)/2 , \\
    \inversekin[{\mathcal{T}}_1]_5 (\theta,x,y,z) & = \! \theta \! - \!
    \inversekin[{\mathcal{T}}_1]_1 (\theta,x,y,z) \! - \!
    \inversekin[{\mathcal{T}}_1]_3 (\theta,x,y,z) ,
  \end{align*}
  where $\rho=\sqrt{\alpha^2 + \beta^2}$ and
  \begin{align*}
    \gamma & = z - d_1 \theta \, ,\\
    \begin{bmatrix}
      \alpha\\
      \beta
    \end{bmatrix}
    &= \frac{1}{b_2^2+c_2^2}
    \begin{bmatrix}
      b_2 & c_2\\
      -c_2 & b_2
    \end{bmatrix}
    \left(
    \begin{bmatrix}
      x \\ y
    \end{bmatrix}
    -
    \begin{bmatrix}
      -c_1 & b_1 \\
      b_1 & c_1
    \end{bmatrix}
    \begin{bmatrix}
      1 - \cos \theta \\
      \sin \theta
    \end{bmatrix}    
  \right).
  \end{align*}
  Then, $\inversekin[{\mathcal{T}}_1]$ is a global right inverse of
  $\forwardkin^{(1,2,1,2,1)}$, that is, it satisfies
  $\map{\forwardkin^{(1,2,1,2,1)} \scirc \inversekin[{\mathcal{T}}_1] =
    \Id_{\SE{2}\times\real}}{\SE{2}\times\real}{\SE{2}\times\real}$.
\end{proposition}

\begin{proof}
  The proof follows from the expression of
  $\forwardkin^{(1,2,1,2,1)}$. Let $\forwardkin^{(1,2,1,2,1)}
  (t_1,t_2,t_3,t_4,t_5) = (\theta,x,y,z)$,
  \begin{align*}
    \theta &= t_1 + t_3 + t_5 \,, \\
    z &= t_2 + t_4 + d_1 \theta \,,\\
    \begin{bmatrix}
      x \\
      y
    \end{bmatrix}
    &=
    \begin{bmatrix}
      -c_1 & b_1 \\
      b_1 & c_1
    \end{bmatrix}
    \begin{bmatrix}
      1 - \cos \theta \\
      \sin \theta
    \end{bmatrix}    
    \\
    & +
    \begin{bmatrix}
      b_2 &- c_2\\
      c_2  & b_2
    \end{bmatrix}
    \left(
      \begin{bmatrix}
        \cos t_1 \\
        \sin t_1
      \end{bmatrix}
      t_2 +
      \begin{bmatrix}
        \cos (t_1 + t_3) \\
        \sin (t_1 + t_3)
      \end{bmatrix}
      t_4 \right) \,.
  \end{align*}
  The equation in $[x,y]^T$ can be rewritten as
  \begin{align*}
    \begin{bmatrix}
      \alpha \\
      \beta
    \end{bmatrix}
    =
    \begin{bmatrix}
      \cos t_1 \\
      \sin t_1
    \end{bmatrix}
    t_2 +
    \begin{bmatrix}
      \cos (t_1+t_3) \\
      \sin (t_1+t_3)
    \end{bmatrix}
    t_4 \, ,
  \end{align*}
  which is solved by the selection of coasting times given by the components
  of the map $\inversekin[{\mathcal{T}}_1]$.
\end{proof}

\begin{proposition} (Inversion for ${\mathcal{T}}_2$-systems on $\SE{2}
  \times \real$).
  \label{prop:SE(2)timesR-case-B}
  Let $(V_1,V_2) \in {\mathcal{T}}_2$. Define the neighborhood of the
  identity in $\SE{2} \times \real$
  \begin{multline*}    
    U = \Big\{ (\theta,x,y,z)\in \SE{2} \times \real\,|\enspace 4 \, \| ( c_1
    -c_2, b_1 - b_2 )\|^2 \geq  \\
    \max\{ \| (x ,y ) \|^2 \,, 2(1 - \cos \theta) \| (b_1,c_1)\|^2 \} \, , \\
    | z - d_1 \theta | \le 2 | d_2 - d_1 | \arccos \left( -1 + \frac{1}{\| (
        c_1 -c_2, b_1 - b_2 )\|} \right. \\ \left.  \cdot \left( \|(x, y)\| +
        \| (b_1,c_1) \| \sqrt{2(1 - \cos \theta)} \right) \right) \Big\}.
  \end{multline*}
  Consider the map $\map{\inversekin[{\mathcal{T}}_2]}{\SE{2} \times
    \real}{\real^5}$ whose components are
  \begin{align*}
    \inversekin[{\mathcal{T}}_2]_1 (\theta,x,y,z) &=
    \atantwo{l}{\sqrt{4 - l^2}} +
    \atantwo{\alpha}{\beta}, \\
    \inversekin[{\mathcal{T}}_2]_2 (\theta,x,y,z) &= 2
    \atantwo{\sqrt{4-l^2}}{l}, \\
    \inversekin[{\mathcal{T}}_2]_3 (\theta,x,y,z) &=
    -\atantwo{\rho-l}{\sqrt{4-(\rho-l)^2}} \\
    & - \inversekin[{\mathcal{T}}_2]_1 (\theta,x,y,z) -
    \inversekin[{\mathcal{T}}_2]_2 (\theta,x,y,z), \\
    \inversekin[{\mathcal{T}}_2]_4 (\theta,x,y,z) &= \gamma -
    \inversekin[{\mathcal{T}}_2]_2 (\theta,x,y,z)  ,\\
    \inversekin[{\mathcal{T}}_2]_5 (\theta,x,y,z) &= \theta -
    \sum_{i=1}^4 \inversekin[{\mathcal{T}}_2]_i (\theta,x,y,z) ,
  \end{align*}
  where $\rho=\sqrt{\alpha^2 + \beta^2}$, $s=\sin (\gamma/2)$, $c=\cos
  (\gamma/2)$ and
  \begin{align*}
    \gamma & = (z - d_1 \theta)/(d_2-d_1) \, ,\\
    l &= \frac{\rho(1+c) + \Sign (\gamma) \sqrt{\rho^2(1 + c)^2 - (1+c)(2
        \rho^2 -8s^2)}}{2(1+c)} , \\
    \begin{bmatrix}
      \alpha\\
      \beta
    \end{bmatrix}
    &= \frac{1}{ \| ( d_1 - d_2, c_1 - c_2 )\|^2}
    \begin{bmatrix}
      d_1 - d_2 & c_2 - c_1\\
      c_1 - c_2 & d_1 - d_2
    \end{bmatrix}
    \\
    & \hspace*{2cm} \cdot \left(
      \begin{bmatrix}
        x \\ y
      \end{bmatrix}
      -
      \begin{bmatrix}
        -d_1 & c_1 \\ 
        c_1 & d_1     
      \end{bmatrix}
      \begin{bmatrix}
        1 - \cos \t \\
        \sin \t
      \end{bmatrix}    
    \right) \, .
  \end{align*}
  Then, $\inversekin[{\mathcal{T}}_2]$ is a local right inverse of
  $\forwardkin^{(1,2,1,2,1)}$, that is, it satisfies
  $\map{\forwardkin^{(1,2,1,2,1)} \scirc \inversekin[{\mathcal{T}}_2] =
    \Id_{U}}{U}{U}$.
\end{proposition}

\begin{proof}
  If $(\theta,x,y,z) \in U$, then $\rho \le 4$ and $|\gamma| \le 2 \arccos
  \left(-1 + \rho/2 \right)$. This in turn implies that
  \begin{align*}
    c = \cos \left( \frac{\gamma}{2} \right) \ge -1 + \frac{\rho}{2} \ge -1 +
    \frac{\rho^2}{8} \,
  \end{align*}
  over $\rho \le 4$. The second inequality guarantees that $l$ is
  well-defined. The first one implies $l \in [\rho-2,2]$, which makes
  $\inversekin[{\mathcal{T}}_2]$ well-defined on $U$.  Let
  $\inversekin[{\mathcal{T}}_2] (\theta,x,y,z)=(t_1,t_2,t_3,t_4,t_5)$.
  The components of $\forwardkin^{(1,2,1,2,1)} (t_1,t_2,t_3,t_4,t_5)$
  are the following
  \begin{align*}
    \theta & = t_1 + t_2 +t_3+t_4 +t_5\,,\\
    z &= d_1 \theta + (d_2 - d_1)
    (t_2 + t_4)\,,\\
      \begin{bmatrix}
        x\\
        y
      \end{bmatrix}
      &=
      \begin{bmatrix}
        -c_1 \\
        b_1
      \end{bmatrix}
      +
      \begin{bmatrix}
        c_1 & b_1 \\
        -b_1 & c_1
      \end{bmatrix}
      \begin{bmatrix}
        \cos \theta\\
        \sin \theta
      \end{bmatrix}
      +
      \begin{bmatrix}
        c_1 - c_2 & b_1 -b_2 \\
        b_2 -b_1 & c_1 -c_2
      \end{bmatrix}
      \\
      & \hspace*{-5pt}
      \begin{bmatrix}
        \cos t_1 - \cos(t_1 + t_2) + \cos(t_1+t_2+t_3) -
        \cos(\sum_{i=1}^4
        t_i)\\
        \sin t_1 - \sin (t_1 + t_2) + \sin (t_1 +t_2 +t_3)
        -\sin(\sum_{i=1}^4 t_i)
      \end{bmatrix}
  \end{align*}
  After some rather involved computations, one can verify 
  $\forwardkin^{(1,2,1,2,1)} (t_1,t_2,t_3,t_4,t_5) = (\theta,x,y,z)$.
\end{proof}

\subsection{Three-dimensional input distribution}

Let $V_1$, $V_2$, $V_3$ satisfy the controllability condition in
Lemma~\ref{lem:ctrl-SE2timesR-three}. Consider
$\map{\forwardkin^{(1,3,2,1)}}{\real^4}{\SE{2} \times \real}$.

\begin{proposition} (Inversion for ${\mathcal{T}}_3$-systems on $\SE{2}
  \times \real$).
  \label{prop:SE(2)timesR-case-C}
  Let $(V_1,V_2,V_3) \in {\mathcal{T}}_3$.  Consider the map
  $\map{\inversekin[{\mathcal{T}}_3]}{\SE{2} \times \real}{\real^4}$ whose
  components are
  \begin{align*}
    \inversekin[{\mathcal{T}}_3]_1 (\theta,x,y,z) &= \atantwo{\alpha}{\beta}
    - \inversekin[{\mathcal{T}}_3]_2
    (\theta,x,y,z)  \,, \\
    \inversekin[{\mathcal{T}}_3]_2 (\theta,x,y,z) &= \frac{z-d_1
      \theta}{d_3-d_1}, \\
    \inversekin[{\mathcal{T}}_3]_3 (\theta,x,y,z) &= \rho  \, , \\
    \inversekin[{\mathcal{T}}_3]_4 (\theta,x,y,z) &= \theta -
    \atantwo{\alpha}{\beta}\, ,
  \end{align*}
  where $\rho=\sqrt{\alpha^2 + \beta^2}$ and
  \begin{align*}
    \begin{bmatrix}
      \alpha\\
      \beta
    \end{bmatrix}
    &= \frac{1}{b_2^2+c_2^2}
    \begin{bmatrix}
      b_2 & c_2\\
      -c_2 & b_2
    \end{bmatrix}
    \left(
    \begin{bmatrix}
      x \\ y
    \end{bmatrix}
    -
    \begin{bmatrix}
      -c_1 & b_1 \\
      b_1 & c_1
    \end{bmatrix}
    \begin{bmatrix}
      1 - \cos \theta \\
      \sin \theta
    \end{bmatrix}    
  \right).
  \end{align*}
  Then, $\inversekin[{\mathcal{T}}_3]$ is a global right inverse of
  $\forwardkin^{(1,3,2,1)}$, that is, it satisfies
  $\map{\forwardkin^{(1,3,2,1)} \scirc \inversekin[{\mathcal{T}}_3] =
    \Id_{\SE{2}\times \real}}{\SE{2}\times \real}{\SE{2}\times \real}$.
\end{proposition}


\begin{proof}
  The proof follows from the expression of $\forwardkin^{(1,3,2,1)}$.
  If $\forwardkin^{(1,3,2,1)} (t_1,t_2,t_3,t_4) = (\theta,x,y,z)$,
  then
  \begin{align*} 
    \theta &= t_1 + t_2 + t_4 \,, \\
    z &= d_1 t_1 + d_3 t_2 + d_1 t_4 = d_1 \theta + (d_3-d_1) t_2 \,, \\
    \begin{bmatrix}
      x \\ y 
    \end{bmatrix}
    & = 
    \begin{bmatrix}
      -c_1 & b_1 \\ 
      b_1 & c_1     
    \end{bmatrix}
    \begin{bmatrix}
      1 - \cos \theta \\
      \sin \theta
    \end{bmatrix}    
    +  
    \begin{bmatrix}
      b_2 & -c_2 \\ 
      c_2 & b_2 
    \end{bmatrix}
    \begin{bmatrix}
      \cos (t_1+t_2) \\ 
      \sin (t_1+t_2)     
    \end{bmatrix} t_3\,.  
  \end{align*} 
  The equation in $[x,y]^T$ can be rewritten as
  \[
  \begin{bmatrix}
    \alpha\\
    \beta
  \end{bmatrix}
  =
  \begin{bmatrix}
    \cos (t_1+t_2)\\
    \sin (t_1+t_2)
  \end{bmatrix}
  t_3 \, ,
  \]
  which is solved by the selection given by $(t_1,t_2,t_3,t_4) =
  \inversekin[{\mathcal{T}}_3](\theta,x,y,z)$.
\end{proof}

Consider the map $\map{\forwardkin^{(1,2,1,3)}}{\real^4}{\SE{2} \times
  \real}$.

\begin{proposition} (Inversion for ${\mathcal{T}}_4$-systems on $\SE{2}
  \times \real$).
  \label{prop:SE(2)-case-D}
  Let $(V_1,V_2,V_3) \in {\mathcal{T}}_4$. Consider the map
  $\map{\inversekin[{\mathcal{T}}_4]}{\SE{2} \times \real}{\real^4}$
  given by
  \begin{multline*}
    \inversekin[{\mathcal{T}}_4] (\theta,x,y,z) = \Big(
    \atantwo{\alpha}{\beta} , \rho, \\ \theta - \atantwo{\alpha}{\beta} ,
    \frac{z-d_1 \theta}{d_3} \Big) ,
  \end{multline*}  
  where $\rho=\sqrt{\alpha^2 + \beta^2}$ and
  \begin{align*}
    \begin{bmatrix}
      \alpha\\
      \beta
    \end{bmatrix}
    &= \frac{1}{b_2^2+c_2^2}
    \begin{bmatrix}
      b_2 & c_2\\
      -c_2 & d_2
    \end{bmatrix}
    \left(
    \begin{bmatrix}
      x \\ y
    \end{bmatrix}
    -
    \begin{bmatrix}
      -c_1 & b_1 \\
      b_1 & c_1
    \end{bmatrix}
    \begin{bmatrix}
      1 - \cos \t \\
      \sin \t
    \end{bmatrix}    
  \right).
  \end{align*}
  Then, $\inversekin[{\mathcal{T}}_4]$ is a global right inverse of
  $\forwardkin^{(1,2,1,3)}$, that is, it satisfies
  $\map{\forwardkin^{(1,2,1,3)} \scirc \inversekin[{\mathcal{T}}_4] =
    \Id_{\SE{2}\times \real}}{\SE{2}\times \real}{\SE{2}\times \real}$.
\end{proposition}


\begin{proof}
  If $\forwardkin^{(1,2,1,3)} (t_1,t_2,t_3,t_4) = (\theta,x,y,z)$, then
  \begin{align*} 
    \theta &= t_1 + t_3 \,, \\
    \begin{bmatrix}
      x \\ y 
    \end{bmatrix}
    & = 
    \begin{bmatrix}
      -c_1 & b_1 \\ 
      b_1 & c_1     
    \end{bmatrix}
    \begin{bmatrix}
      1 - \cos \t \\
      \sin \t
    \end{bmatrix}    
    +  
    \begin{bmatrix}
      b_2 & -c_2 \\ 
      c_2 & b_2 
    \end{bmatrix}
    \begin{bmatrix}
      \cos t_1 \\ 
      \sin t_1     
    \end{bmatrix} t_2 \,, \\
    z & = d_1 (t_1 + t_3) + d_3 t_4 \, .
  \end{align*} 
  The equation in $[x,y]^T$ can be rewritten as $[\alpha, \beta]^T = [ \cos
  t_1, \sin t_1]^T t_2$. As in the proof of
  Proposition~\ref{prop:SE(2)-case-A}, the selection $t_1
  =\atantwo{\alpha}{\beta}$, $t_2 = \rho$ solves it.
\end{proof}

\begin{proposition} (Inversion for ${\mathcal{T}}_5$-systems on $\SE{2}
  \times \real$).
  \label{prop:SE(2)timesR-case-E}
  Let $(V_1,V_2,V_3) \in {\mathcal{T}}_5$. Define the neighborhood of
  the identity in $\SE{2} \times \real$
  \begin{multline*}    
    U = \{(\theta,x,y)\in \SE{2} \times \real\,|\enspace \| ( c_1 - c_2, b_1
    - b_2 )\|^2\geq
    \\
    \max\{ \| (x ,y ) \|^2 \,, 2(1 - \cos \theta) \| (b_1,c_1)\|^2 \}.
  \end{multline*}
  Consider the map $\map{\inversekin[{\mathcal{T}}_5]}{U\subset\SE{2} \times
    \real}{\real^4}$ whose components are
  \begin{align*}
    \inversekin[{\mathcal{T}}_5]_1 (\theta,x,y,z) &= \atantwo{\rho}{\sqrt{4
        -\rho^2}} +
    \atantwo{\alpha}{\beta} \,, \\
    \inversekin[{\mathcal{T}}_5]_2 (\theta,x,y,z) &= \atantwo{2-\rho^2}{\rho
      \sqrt{4-\rho^2}}
    \!, \\
    \inversekin[{\mathcal{T}}_5]_3 (\theta,x,y,z) &= \theta
    -\inversekin[{\mathcal{T}}_5]_1
    (\theta,x,y)-\inversekin[{\mathcal{T}}_5]_2 (\theta,x,y) \, , \\
    \inversekin[{\mathcal{T}}_5]_4 (\theta,x,y,z) &= \frac{z-d_1 \theta}{d_3}
    \, ,
  \end{align*}
  and $\rho=\sqrt{\alpha^2 + \beta^2}$ and
  \begin{align*}
    \begin{bmatrix}
      \alpha\\
      \beta
    \end{bmatrix}
    &= \frac{1}{ \| ( c_1 - c_2, b_1 - b_2 )\|^2}
    \begin{bmatrix}
      c_1 - c_2 & b_2 - b_1\\
      b_1 - b_2 & c_1 - c_2
    \end{bmatrix}
    \\
    & \hspace*{2cm} \cdot \left(
    \begin{bmatrix}
      x \\ y
    \end{bmatrix}
    -
    \begin{bmatrix}
      -c_1 & b_1 \\
      b_1 & c_1
    \end{bmatrix}
    \begin{bmatrix}
      1 - \cos \theta \\
      \sin \theta
    \end{bmatrix}    
  \right) \, .
  \end{align*}
  Then, $\inversekin[{\mathcal{T}}_5]$ is a local right inverse of
  $\forwardkin^{(1,2,1,3)}$, that is, it satisfies
  $\map{\forwardkin^{(1,2,1,3)} \scirc \inversekin[{\mathcal{T}}_5] =
    \Id_{U}}{U}{U}$.
\end{proposition}

\begin{proof}
  If $(\theta,x,y,z) \in U$, then one can see that $\rho = \|
  (\alpha,\beta)\| \leq 2$, and therefore $\inversekin[{\mathcal{T}}_5]$ is
  well-defined on~$U$.  Let $\inversekin[{\mathcal{T}}_5]
  (\theta,x,y,z)=(t_1,t_2,t_3,t_4)$. The components of
  $\forwardkin^{(1,2,1,3)} (t_1,t_2,t_3,t_4)$ are
  \begin{align*}
    & \forwardkin^{(1,2,1,3)}_1 (t_1,t_2,t_3,t_4) = t_1 + t_2 + t_3 \,, \\
    & \begin{bmatrix}
      \forwardkin^{(1,2,1,3)}_2 (t_1,t_2,t_3,t_4) \\
      \forwardkin^{(1,2,1,3)}_3 (t_1,t_2,t_3,t_4)
    \end{bmatrix}
    =
    \begin{bmatrix}
      -c_1 & b_1 \\
      b_1 & c_1
    \end{bmatrix}
    \begin{bmatrix}
      1 - \cos \theta  \\
      \sin \theta
    \end{bmatrix}
    \\
    & \hspace*{1.5cm} +
    \begin{bmatrix}
      c_1-c_2 &  b_1-b_2  \\
      b_2- b_1 & c_1 - c_2
    \end{bmatrix}
    \begin{bmatrix}
      \cos t_1 - \cos (t_1 + t_2)\\
      \sin t_1 - \sin (t_1 + t_2)
    \end{bmatrix} \,,\\
    & \forwardkin^{(1,2,1,3)}_4 (t_1,t_2,t_3,t_4) = d_1(t_1 + t_2 + t_3) +
    d_3 t_4 \,.
  \end{align*}
  One can verify that $\forwardkin^{(1,2,1,3)} (t_1,t_2,t_3,t_4) =
  (\theta,x,y,z)$.
\end{proof}

\begin{figure}[htbp]
  \centering
  \includegraphics[width=.9\linewidth]{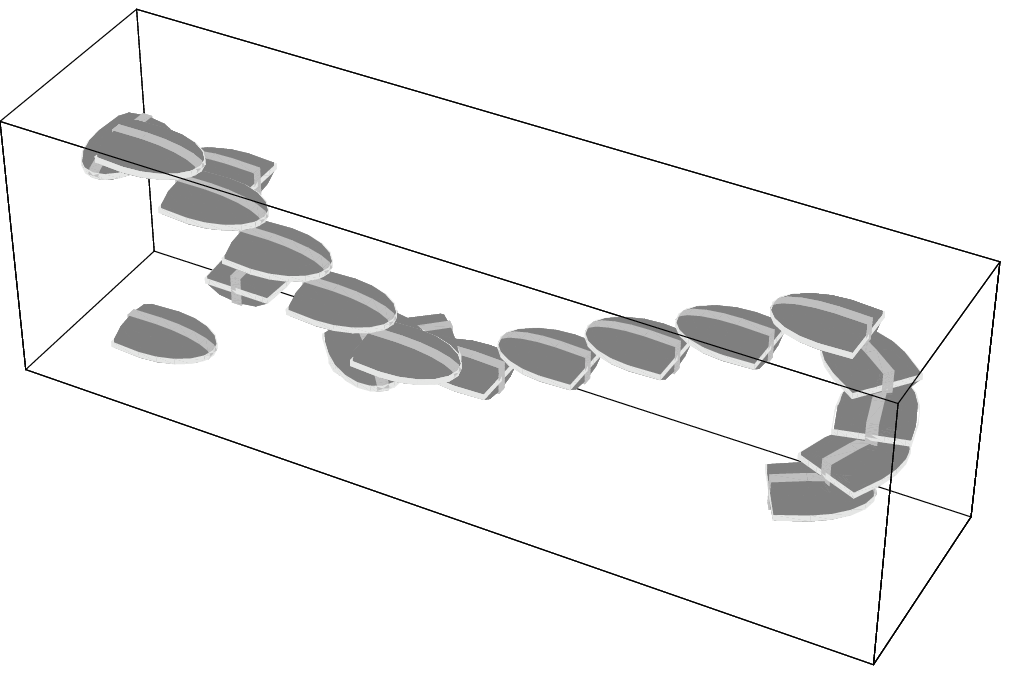}
  \caption{We illustrate the inverse-kinematics planner for 
    a ${\mathcal{T}}_1$-system on $\SE{2}\times\real$.  The system
    parameters are $b_1=1$, $c_1=0$, $d_1=.5$, $b_2=-2$, and $c_2=0$.
    The target location is $(\pi/6,10,0,1)$.}
  \label{fig:se2real}
\end{figure}

We illustrate the performance of the algorithms in Fig.~\ref{fig:se2real}.

\section{Conclusions}
We have presented a catalog of feasible motion planning algorithms for
underactuated controllable systems on $\SE{2}$, $\SO{3}$ and
$\SE{2}\times\real$. Future directions of research include (i) considering
other relevant classes of underactuated systems on $\SE{3}$, (ii) computing
catalogs of optimal sequences of motion primitives, and (iii) developing
hybrid feedback schemes that rely on the proposed open-loop planners to
achieve point stabilization and trajectory tracking.

\subsection*{Acknowledgments}
The authors would like to thank Professor Kevin M. Lynch for insightful
discussions.  This research was partially supported by NSF Grant
IIS-0118146 and Spanish MCYT Grant BFM2001-2272.

\end{document}